\documentclass[10pt]{article}

\usepackage{amssymb,amsmath,latexsym}
\usepackage[dvips]{graphics}
\usepackage{mathrsfs}
\usepackage{extarrows}

\newtheorem{theorem}{Theorem}[section]
\newtheorem{prop}[theorem]{Proposition}

\newtheorem{cor}[theorem]{Corollary}

\newtheorem{lemma}[theorem]{Lemma}
\numberwithin{equation}{section}

\def\pf{\noindent{\bf Proof.} }
\def\qed{{\hfill $\Box$\medskip}}

\def\PP{\mathbb{P}}
\def\EE{\mathbb{E}}

\def\dd{\mathrm{d}}
\def\Beta{\mathrm{B}}
\def\id{\mathbf{1}}
\def\one{\uppercase\expandafter{\romannumeral1}}
\def\two{\uppercase\expandafter{\romannumeral2}}
\allowdisplaybreaks[4]

\input{epsf.sty}
\usepackage{epsfig}

\begin{document}
\title{\bf\Large A new method for computing the expected hitting time between arbitrary different configurations of the multiple--urn Ehrenfest model}
\author{Sai SONG\footnote{Key Laboratory of Advanced Theory and Application in Statistics and Data Science-MOE, School of Statistics, East China Normal University.}~~and Qiang YAO\footnote{Corresponding author. Key Laboratory of Advanced Theory and Application in Statistics and Data Science-MOE, School of Statistics, East China Normal University and NYU--ECNU Institute of Mathematical Sciences at NYU Shanghai.
                      E-mail: qyao@sfs.ecnu.edu.cn.}
                      }

\maketitle{}

\begin{abstract}
We study a multiple--urn version of the Ehrenfest model. In this setting, we denote the $n$ urns by Urn $1$ to Urn $n$, where $n\geq2$. Initially, $M$ balls are randomly placed in the $n$ urns. At each subsequent step, a ball is selected and put into the other $n-1$ urns with equal probability. The expected hitting time leading to a change of the $M$ balls' status is computed using the method of stopping times. As a corollary, we obtain the expected hitting time of moving all the $M$ balls from Urn $1$ to Urn $2$. This proves a conjecture which was recently made in Chen et al.~\cite{Chen-Goldstein-Lathrop-Nelsen2017}.
\end{abstract}
\noindent{\bf 2010 MR subject classification:}~60C05, 60J10


\section{Introduction}

We extend the classical two--urn Ehrenfest model to the multiple--urn case. Label the $n$ urns by Urn $1$ to Urn $n$, where $n\geq2$. At the beginning, $M$ balls are arbitrarily placed in the $n$ urns. Then at each time, one ball is chosen at random, taken from the current urn it resides in, and placed in one of the other $n-1$ urns with equal probability. This model can be treated as a symmetric simple random walk on the graph $G_M=(V_M,E_M)$, where $V_M=\{1,\ldots,n\}^M$, and $E$ contains edges connecting two vertices in $V_M$ if exactly one of their components differs. Here the subscript ``$M$'' is to stress that the number of balls is $M$. Therefore, $G_M$ is a transitive graph~(that is, for any $e,e^\prime\in E_M$, there is an automorphism of the graph that takes $e$ to $e^\prime$) with $n^M$ vertices, and each vertex has common degree $(n-1)M$. Strictly speaking, if we let $X_t=(X_t^{(1)},\ldots,X_t^{(M)})$ be the state at time $t=0,1,\ldots$, where $X_t^{(i)}$ is the number of the urn in which the $i$th ball resides at time $t$, then $\{X_t:~t=0,1,\ldots\}$ is a time homogeneous Markov chain on $V_M$ with transition probability
\begin{align}\label{e:transitionRW}
&p_{(x_1,\ldots,x_M),(y_1,\ldots,y_M)}\nonumber\\
=&\left\{\begin{array}{ll}\frac{1}{(n-1)M},~~~~~~~~~~~~\text{if there exists}~i~\text{s.t.}~x_i\neq y_i,~\text{and}~x_j=y_j~\text{for}~j\neq i;\\~~~~~0,~~~~~~~~~~~~~~~~~~~~~~~~~~~~~~~~~~~~~~~~~\text{otherwise}.\end{array}\right.
\end{align}
For $x_1,\ldots,x_M\in\{1,2,\ldots,n\}$, denote by $$T_{(x_1,\ldots,x_M)}=\inf\{t\geq0:~X_t=(x_1,\ldots,x_M)\}$$ the first time that $\{X_t\}$ hits state $(x_1,\ldots,x_M)$. Our main result is described in the following theorem.
\begin{theorem}\label{t:main2}
For any two different configurations $(a_1,\ldots,a_M),(b_1,\ldots,b_M)\in\{1,\ldots,n\}^M$, denote $L=\sum\limits_{i=1}^M\id_{\{a_i=b_i\}}$. Then $$\EE(T_{(b_1,\ldots,b_M)}~|~X_0=(a_1,\ldots,a_M))=\sum\limits_{k=L}^{M-1}\frac{(n-1)^{k+1}}{{M-1 \choose k}}\sum\limits_{i=0}^k\frac{{M \choose i}}{(n-1)^i},$$ where $\id_{\{\cdot\}}$ denotes the indicator function, and ${n \choose m}:=\dfrac{n!}{m!(n-m)!}~~(0\leq m\leq n)$ denotes the combinatorial number.
\end{theorem}

\medskip

Note that when $L=M-1$, the righthand side of Theorem \ref{t:main2} becomes $\displaystyle(n-1)^M\sum\limits_{i=0}^{M-1}\frac{{M \choose i}}{(n-1)^i}$, which equals to $n^M-1$. This is a well known result in Markov chains.

\bigskip

As a special case, we obtain the following corollary, which provides the expected hitting time of moving all balls from Urn $1$ to Urn $2$.
\begin{cor}\label{c:main}
$\EE\left(\left.T_{(\tiny{\underbrace{2,2,\ldots,2}_M})}~\right|~X_0=(\footnotesize{\underbrace{1,1,\ldots,1}_M})\right)=\dfrac{(n-1)M}{n}\sum\limits_{k=1}^M\dfrac{n^k}{k}$.
\end{cor}

\bigskip

\noindent\textbf{Remark.}~(1)~Chen et al.~\cite{Chen-Goldstein-Lathrop-Nelsen2017} proved Corollary \ref{c:main} for the special case $n=3$ by using the method of electric networks. They conjectured that the result for general multiple--urn case should be of the form as stated in Corollary \ref{c:main}.

(2)~Corollary \ref{c:main} is a special case of Theorem \ref{t:main2} by letting $L=0$. This is not a straightforward result. A key step is to establish the equality
\begin{equation}\label{e:main}
\sum\limits_{k=0}^{M-1}\frac{(n-1)^k}{M{M-1 \choose k}}\sum\limits_{i=0}^k\frac{{M \choose i}}{(n-1)^i}=\frac{1}{n}\sum\limits_{k=1}^M\frac{n^k}{k}.
\end{equation}
It is true. Although it is well known that one cannot get an explicit formula for the inner sum $\sum\limits_{i=0}^k\dfrac{{M \choose i}}{(n-1)^i}$ for fixed $k$, we can simplify the double sum on the left side of (\ref{e:main}) to the form on the right side. In fact, by direct calculation, we obtain
\begin{align*}
&\text{LHS of}~(\ref{e:main})=\sum\limits_{k=0}^{M-1}\sum\limits_{i=0}^k(n-1)^{k-i}\frac{{M \choose i}}{M{M-1 \choose k}}\\
=&\sum\limits_{k=0}^{M-1}\sum\limits_{i=0}^k{k \choose i}(n-1)^{k-i}\cdot\Beta(k-i+1,M-k)~~~~~~(\Beta(\cdot,\cdot)~\text{denotes the Beta function})\\
=&\sum\limits_{k=0}^{M-1}\sum\limits_{i=0}^k{k \choose i}(n-1)^{k-i}\cdot\int_0^1s^{k-i}(1-s)^{M-k-1}\dd s\\
=&\sum\limits_{k=0}^{M-1}\int_0^1[1+(n-1)s]^k(1-s)^{M-k-1}\dd s\\
=&\int_0^1\frac{[1+(n-1)s]^M-(1-s)^M}{ns}\dd s\\
=&\sum\limits_{k=1}^M\int_0^1{M \choose k}(ns)^{k-1}(1-s)^{M-k}\dd s\\
=&\sum\limits_{k=1}^M{M \choose k}n^{k-1}\cdot\Beta(k,M-k+1)=\sum\limits_{k=1}^M\frac{n^{k-1}}{k}=\text{RHS of}~(\ref{e:main}).
\end{align*}
Therefore, (\ref{e:main}) is obtained.

\bigskip

The Ehrenfest model was proposed in Ehrenfest \& Ehrenfest~\cite{Ehrenfest-Ehrenfest1907} for the first time. It is ``a test bed of key concepts of statistical mechanics''~(see Meerson \& Zilber~\cite{Meerson-Zilber2018}). There are many problems concerning this simple but insightful model. The study of the hitting time was first restricted in the $2$--urn case, see Blom~\cite{Blom1989}, Lathrop et al.~\cite{Lathrop-Goldstein-Chen2016}, Palacios~\cite{Palacios1994}, etc. The multiple--urn model has attracted much attention recently in both theoretical and application fields, see for example, Xue~\cite{Xue2021}, Aloisi \& Nali~\cite{Aloisi-Nali2018}, etc. So it is worthwhile to investigate the model more deeply.

\medskip

This paper is devoted to a purely probabilistic proof for Theorem \ref{t:main2}. The method we used is original as far as we know. There are other alternative methods available. For example, the method of electric networks can be used. Readers can refer to Doyle and Snell~\cite{Doyle-Snell1984} or Lyons and Peres~\cite{Lyons-Peres2017} for the introduction of this method. Palacios~\cite{Palacios1994} used the method of electric networks to consider the $2$--urn case. Chen et al.~\cite{Chen-Goldstein-Lathrop-Nelsen2017} used this method to consider the $3$--urn case. However, it may become difficult to use the ``Y--Delta'' transformation to simplify the network when the number of urns $n$ is large, especially when $n\geq4$. We mention here that Chen et al.~\cite{Chen-Goldstein-Lathrop-Nelsen2017} stated the main difficulty in the conclusion section of their paper. Besides, the method of the auxiliary Markov chain~(which only considers the number of balls in one fixed urn) inspired by Blom~\cite{Blom1989} and Lathrop et al.~\cite{Lathrop-Goldstein-Chen2016} can also settle this problem. Very recently, Xin et al.~\cite{Xin-Zhao-Yao-Cui2020} obtained the distribution of the hitting time, which can imply its expectation, too. The strength of our present method is to make a good illustration for the recursive formulas (\ref{e:sk_induction}) and (\ref{e:sk0_induction}) below, especially for the special case (\ref{e:sk0_induction}).

\bigskip

The rest of this paper is organized as follows. In Sections 2, we will utilize some stopping times to set up the recursive relation. In Section 3, we will use the above recursive relationship to prove Theorem \ref{t:main2}. And in Section 4, we will give some concluding remarks.

\section{The recursive relation}

We fix the number of urns~($n$) and let the number of balls~($k$) vary. If we consider $k$ balls, then $\{X_t\}$ becomes a symmetric simple random walk on the graph $G_k=(V_k,E_k)$, that is, the state space of $\{X_t\}$ becomes $V_k=\{1,\ldots,n\}^k$, and the transition probability is
\begin{align*}
&p_{(x_1,\ldots,x_k),(y_1,\ldots,y_k)}\nonumber\\
=&\left\{\begin{array}{ll}\frac{1}{(n-1)k},~~~~~~~~~~~~\text{if there exists}~i~\text{s.t.}~x_i\neq y_i,~\text{and}~x_j=y_j~\text{for}~j\neq i;\\~~~~0,~~~~~~~~~~~~~~~~~~~~~~~~~~~~~~~~~~~~~~~~~\text{otherwise}.\end{array}\right.
\end{align*}
For any $L\geq0$, define
$$s_{L,k}:=\EE\left(\left.T_{(\tiny{\underbrace{2,2,\ldots,2}_k})}~\right|~X_0=(\footnotesize{\underbrace{2,\ldots,2}_L,\underbrace{1,\ldots,1}_{k-L}})\right)$$
for $k\geq L$. By the transitivity of $G_k$, if $\sum\limits_{i=1}^k\id_{\{a_i=b_i\}}=L$, then
$$\EE(T_{(b_1,\ldots,b_k)}~|~X_0=(a_1,\ldots,a_k))=s_{L,k}.$$

\bigskip

Our main step is to prove the following proposition, which gives an recursive formula for $\{s_{L,k}:~k=L,L+1,\ldots\}$, where $L\geq0$ is fixed at first.
\begin{prop}\label{p:sk_induction}
$\{s_{L,k}:~k=L,L+1,\ldots\}$ satisfies the recursive formula
\begin{equation}\label{e:sk_induction}
\begin{cases}
\displaystyle s_{L,k}=\frac{k}{k-1}s_{L,k-1}+(n-1)n^{k-1}+\frac{n-1}{{k-1 \choose L}}\sum\limits_{i=1}^L{k-1 \choose i-1}(n-1)^{L-i}, & k\geq L+1;\\
s_{L,L}=0. &
\end{cases}
\end{equation}
\end{prop}

\noindent\textbf{Remark.}~The intuitive description for the recursive formula (\ref{e:sk_induction}) is as follows. If $\{X_t\}$ hits $(\footnotesize{\underbrace{2,2,\ldots,2}_k})$ from $(\footnotesize{\underbrace{2,\ldots,2}_L,\underbrace{1,\ldots,1}_{k-L}})$, the first $k-1$ components~(denoted by $\{Y_t\}$, whose rigorous definition will be given in the proof of Lemma \ref{l:Etau_1}) must hit $(\footnotesize{\underbrace{2,2,\ldots,2}_{k-1}})$ from $(\footnotesize{\underbrace{2,\ldots,2}_L,\underbrace{1,\ldots,1}_{k-L-1}})$ first. Since $\{Y_t\}$ can be seen as a ``delayed'' random walk on $V_{k-1}$, the above procedure will take $\frac{k}{k-1}s_{L,k-1}$ steps on average. Then it will take some more steps to reach the destination.

\bigskip

The rest of this section is devoted to a strict proof of (\ref{e:sk_induction}). For any $(x_1,\ldots,x_{k-1})\in V_{k-1}$, define
$$A_{(x_1,\ldots,x_{k-1})}:=\{(x_1,\ldots,x_{k-1},x_k):~x_k\in\{1,\ldots,n\}\}.$$
That is, $A_{(x_1,\ldots,x_{k-1})}$ contains all points in $V_k$ such that the first $k-1$ components are $x_1,\ldots,x_{k-1}$. Now we define several stopping times. First, define
$$T_k=\inf\left\{t>0:~X_t=(\footnotesize{\underbrace{2,2,\ldots,2}_k})\right\}.$$ Then
\begin{equation}\label{e:sk_defmethod2}
s_{L,k}=\EE_{(\tiny{\underbrace{2,\ldots,2}_L,\underbrace{1,\ldots,1}_{k-L}})}(T).
\end{equation}
Next, for any $k>L$, we define a sequence of stopping times $\{\tau_{k,i}:~k=0,1,2,\ldots\}$ inductively by $\tau_{k,0}=0$ and
$$\tau_{k,i}=\inf\left\{t>\tau_{k,i-1}:~X_t\in A_{(\tiny{\underbrace{2,\ldots,2}_{k-1}})}\right\}$$
for $i\geq1$. Clearly, for any $k>L$,
$$\PP_{(\tiny{\underbrace{2,\ldots,2}_L,\underbrace{1,\ldots,1}_{k-L}})}\left(\bigcup\limits_{i=1}^\infty\{T_k=\tau_{k,i}\}\right)=\sum\limits_{i=1}^\infty \PP_{(\tiny{\underbrace{2,\ldots,2}_L,\underbrace{1,\ldots,1}_{k-L}})}(T_k=\tau_{k,i})=1.$$\\

The next two lemmas give some important properties of the above stopping times.
\begin{lemma}\label{l:Etau_1}
$\displaystyle\EE_{(\tiny{\underbrace{2,\ldots,2}_L,\underbrace{1,\ldots,1}_{k-L}})}(\tau_{k,i}-\tau_{k,i-1})=\left\{\begin{array}{ll}\frac{k}{k-1}s_{L,k-1},~~~~~~~~~\text{if}~~i=1;\\~~~~n^{k-1},~~~~~~~~~~~\text{if}~~i\geq2.\end{array}\right.$
\end{lemma}
\pf We define two auxiliary Markov chains on $V_{k-1}$. First, denote
$$Y_t:=(X_t^{(1)},\ldots,X_t^{(k-1)})$$ for $t=0,1,2,\ldots$. Then $\{Y_t\}$ is a Markov chain on $V_{k-1}$ which illustrates the positions of the first $k-1$ balls with transition probability
$$p_{xy}=\left\{\begin{array}{ll}~~~~\frac{1}{k},~~~~~~~~~~~\text{if}~y=x;\\~\frac{1}{k(n-1)},~~~~~~~~\text{if}~y\sim x;\\~~~~0,~~~~~~~~~~~\text{otherwise.} \end{array}\right.$$
Here the notation ``$y\sim x$'' means that $y$ is a neighbor of $x$.

Next, denote by $\{Z_t\}$ the Markov chain on $V_{k-1}$ which illustrates the process with $k-1$ balls and $n$ urns with transition probability
$$q_{xy}=\left\{\begin{array}{ll}~\frac{1}{(k-1)(n-1)},~~~~~~~~\text{if}~y\sim x;\\~~~~~~0,~~~~~~~~~~~~~~\text{otherwise.} \end{array}\right.$$

For $x\in V_{k-1}$, let $f_x=\EE_x(T_{k-1})$ for $\{Y_t\}$~(that is, under the transition probability $\{p_{xy}\}$). Let $g_x=\EE_x(T_{k-1})$ for $\{Z_t\}$~(that is, under the transition probability $\{q_{xy}\}$). By the transition probability $\{q_{xy}\}$, we know that $\{g_x\}$ satisfies
\begin{equation}\label{e:g_induction}
g_x=\left\{\begin{array}{ll}1+\frac{1}{(k-1)(n-1)}\sum\limits_{y\sim x}g_y,~~~~~~~~~\text{if}~x\neq(\footnotesize{\underbrace{2,\ldots,2}_{k-1}});\\~~~~~~~~~~~~0,~~~~~~~~~~~~~~~~~~~~~\text{if}~x=(\footnotesize{\underbrace{2,\ldots,2}_{k-1}}). \end{array}\right.
\end{equation}
Similarly, by the transition probability $\{p_{xy}\}$, $\{f_x\}$ satisfies
\begin{equation}\label{e:f_induction}
f_x=\left\{\begin{array}{ll}1+\frac{1}{k}\cdot f_x+\frac{1}{k(n-1)}\sum\limits_{y\sim x}f_y,~~~~~~~~~\text{if}~x\neq(\footnotesize{\underbrace{2,\ldots,2}_{k-1}});\\~~~~~~~~~~~~~~0,~~~~~~~~~~~~~~~~~~~~~~~~~\text{if}~x=(\footnotesize{\underbrace{2,\ldots,2}_{k-1}}). \end{array}\right.
\end{equation}
Note that (\ref{e:f_induction}) can be written as
\begin{equation}\label{e:f_induction1}
\frac{k-1}{k}f_x=\left\{\begin{array}{ll}1+\frac{1}{(k-1)(n-1)}\sum\limits_{y\sim x}\frac{k-1}{k}f_y,~~~~~~~~~~\text{if}~x\neq(\footnotesize{\underbrace{2,\ldots,2}_{k-1}});\\~~~~~~~~~~~~~~~0,~~~~~~~~~~~~~~~~~~~~~~~\text{if}~x=(\footnotesize{\underbrace{2,\ldots,2}_{k-1}}). \end{array}\right.
\end{equation}
Comparing (\ref{e:f_induction1}) with (\ref{e:g_induction}), we can see that $\{g_x\}$ and $\left\{\dfrac{k-1}{k}f_x\right\}$ obey the same difference equation and have the same initial value. Therefore, $g_x=\dfrac{k-1}{k}f_x$ for any $x\in V_{k-1}$. Especially, $f_{(\tiny{\underbrace{2,\ldots,2}_L,\underbrace{1,\ldots,1}_{k-L-1}})}=\dfrac{k}{k-1}g_{(\tiny{\underbrace{2,\ldots,2}_L,\underbrace{1,\ldots,1}_{k-L-1}})}$.

Since $\EE_{(\tiny{\underbrace{2,\ldots,2}_L,\underbrace{1,\ldots,1}_{k-L-1}},a)}(\tau_{k,1})$ is the same for any $a\in\{1,\ldots,n\}$, we obtain $$f_{(\tiny{\underbrace{2,\ldots,2}_L,\underbrace{1,\ldots,1}_{k-L-1}})}=\EE_{(\tiny{\underbrace{2,\ldots,2}_L,\underbrace{1,\ldots,1}_{k-L}})}(\tau_{k,1}).$$ Together with the fact that $g_{(\tiny{\underbrace{2,\ldots,2}_L,\underbrace{1,\ldots,1}_{k-L-1}})}=s_{L,k-1}$, we get $\EE_{(\tiny{\underbrace{2,\ldots,2}_L,\underbrace{1,\ldots,1}_{k-L}})}(\tau_{k,1})=\frac{k}{k-1}s_{L,k-1}$.

\bigskip

Furthermore, since $\{Y_t\}$ is a reversible Markov chain on $V_{k-1}$~(with $n^{k-1}$ vertices), there exists a unique invariant distribution which puts an equal mass of $\dfrac{1}{n^{k-1}}$ on the $n^{k-1}$ vertices. So $\EE_{(\tiny{\underbrace{2,\ldots,2}_L,\underbrace{1,\ldots,1}_{k-L}})}(\tau_{k,i}-\tau_{k,i-1})=n^{k-1}$ for any $i=2,3,\ldots$, as desired.\qed\\

\begin{lemma}\label{l:same_dist}
For any $i=1,2,\ldots$ and $t\geq0$,
$$\PP_{(\tiny{\underbrace{2,\ldots,2}_L,\underbrace{1,\ldots,1}_{k-L}})}(T_k-\tau_{k,i+1}=t~|~T_k>\tau_{k,i})=\PP_{(\tiny{\underbrace{2,\ldots,2}_{k-1},1})}(T_k-\tau_{k,1}=t).$$
That is, $\{(T_k-\tau_{k,i+1}~|~T_k>\tau_{k,i})\}$ have the same distribution for $i=1,2,\ldots$.
\end{lemma}
\pf Note that $\{T_k>\tau_{k,i}\}=\displaystyle\bigcup\limits_{x\neq2}\left\{X_{\tau_{k,i}}=(\footnotesize{\underbrace{2,\ldots,2}_{k-1},x})\right\}$ for any $i=1,2,\ldots$. From the strong Markov property, for any $t\geq0$, $i=1,2,\ldots$ and $x\neq2$,
\begin{align*}
&\PP_{(\tiny{\underbrace{2,\ldots,2}_L,\underbrace{1,\ldots,1}_{k-L}})}\left(T_k-\tau_{k,i+1}=t~\left|~X_{\tau_{k,i}}=(\footnotesize{\underbrace{2,\ldots,2}_{k-1},x})\right.\right)\\
=&\sum\limits_{m=1}^\infty \PP_{(\tiny{\underbrace{2,\ldots,2}_L,\underbrace{1,\ldots,1}_{k-L}})}\left(T_k-\tau_{k,i}=m+t,~\tau_{k,i+1}-\tau_{k,i}=m~\left|~X_{\tau_{k,i}}=(\footnotesize{\underbrace{2,\ldots,2}_{k-1},x})\right.\right)\\
=&\sum\limits_{m=1}^\infty \PP_{(\tiny{\underbrace{2,\ldots,2}_{k-1},x})}(T_k=m+t,~\tau_{k,1}=m)\\
=&\PP_{(\tiny{\underbrace{2,\ldots,2}_{k-1},x})}(T_k-\tau_{k,1}=t)=\PP_{(\tiny{\underbrace{2,\ldots,2}_{k-1},1})}(T_k-\tau_{k,1}=t).
\end{align*}
The last equality is due to the transitivity of $G_k$. Therefore, $$\PP_{(\tiny{\underbrace{2,\ldots,2}_L,\underbrace{1,\ldots,1}_{k-L}})}(T_k-\tau_{k,i+1}=t~|~T_k>\tau_{k,i})=\PP_{(\tiny{\underbrace{2,\ldots,2}_{k-1},1})}(T_k-\tau_{k,1}=t)$$ for any $i=1,2,\ldots$ and $t\geq0$, as desired.\qed\\

For $i=1,2,\ldots,k$, denote
$$B_{2i-1}=\left\{(x_1,\ldots,x_k)\in V_k:~\sum\limits_{\ell=1}^{k-1}\id_{\{x_\ell=2\}}=i-1,~x_k\neq2\right\},$$
$$B_{2i}=\left\{(x_1,\ldots,x_k)\in V_k:~\sum\limits_{\ell=1}^{k-1}\id_{\{x_\ell=2\}}=i-1,~x_k=2\right\},$$
where $\id_{\{\cdot\}}$ is an indicator function. Then $V_k=\displaystyle\bigcup\limits_{m=1}^{2k}B_m$. By the transitivity of $G_k$, $\PP_x(T_k=\tau_{k,1})$ is the same for the $x$'s belonging to the same $B_m$. So $\PP_x(T_k=\tau_{k,1})=\PP_{B_m}(T_k=\tau_{k,1})$ for any $m=1,\ldots,2k$ and any $x\in B_m$. Denote $p_m=\PP_{B_m}(T_k=\tau_{k,1})$ for $m=1,\ldots,2k$. Note that we do not need to emphasize $k$ in the definition of $B_m$'s and $p_n$'s since $k$ is fixed now.

\bigskip

Having defined the $p_n$'s, we can write down the ``embryonic'' version of the recursive formula (\ref{e:sk_induction}).
\begin{lemma}\label{l:sk_induction2}
For any $L\geq0$, $\{s_{L,k}:~k=L,L+1,\ldots\}$ satisfies
\begin{equation}\label{e:sk_induction2}
\begin{cases}
\displaystyle s_{L,k}=\frac{k}{k-1}s_{L,k-1}+\frac{1-p_{2L+1}}{1-p_{2k}}\cdot(n-1)n^{k-1}, & k\geq L+1;\\
s_{L,L}=0. &
\end{cases}
\end{equation}
\end{lemma}

\pf $s_{L,L}=0$ is obvious. When $k\geq L+1$, let
$$u=\EE_{(\tiny{\underbrace{2,\ldots,2}_L,\underbrace{1,\ldots,1}_{k-L}})}(T_k-\tau_{k,1}),~~~~v=\EE_{(\tiny{\underbrace{2,\ldots,2}_L,\underbrace{1,\ldots,1}_{k-L}})}(T_k-\tau_{k,2}~|~T_k>\tau_{k,1}).$$ Then by Lemma \ref{l:same_dist}, $\EE_{(\tiny{\underbrace{2,\ldots,2}_L,\underbrace{1,\ldots,1}_{k-L}})}(T_k-\tau_{k,3}~|~T_k>\tau_{k,2})=v$. Similar to the proof of Lemma \ref{l:same_dist}, we can get
\begin{align*}
&\EE_{(\tiny{\underbrace{2,\ldots,2}_L,\underbrace{1,\ldots,1}_{k-L}})}\left(T_k-\tau_{k,1}~\left|~X_{\tau_{k,1}}=(\footnotesize{\underbrace{2,\ldots,2}_{k-1},x})\right.\right)\\
=&\left\{\begin{array}{ll}\EE_{(\tiny{\underbrace{2,\ldots,2}_L,\underbrace{1,\ldots,1}_{k-L}})}(T_k-\tau_{k,1}~|~T_k>\tau_{k,1}),~~~~~~~~~~\text{if}~x\neq2;\\~~~~~~~~~~~~~~~~~~~~~~0,~~~~~~~~~~~~~~~~~~~~~~~~~~~~~~~\text{if}~x=2.\end{array}\right.
\end{align*}
Therefore,
\begin{align}\label{e:pre1}
&\EE_{(\tiny{\underbrace{2,\ldots,2}_L,\underbrace{1,\ldots,1}_{k-L}})}(T_k-\tau_{k,1})\nonumber\\
=&\sum\limits_{x\neq2}\EE_{(\tiny{\underbrace{2,\ldots,2}_L,\underbrace{1,\ldots,1}_{k-L}})}\left(T_k-\tau_{k,1}~\left|~X_{\tau_{k,1}}=(\footnotesize{\underbrace{2,\ldots,2}_{k-1},x})\right.\right)\cdot \PP\left(X_{\tau_{k,1}}=(\footnotesize{\underbrace{2,\ldots,2}_{k-1},x})\right)\nonumber\\
=&\EE_{(\tiny{\underbrace{2,\ldots,2}_L,\underbrace{1,\ldots,1}_{k-L}})}(T_k-\tau_{k,1}~|~T_k>\tau_{k,1})\cdot\sum\limits_{x\neq2}\PP\left(X_{\tau_{k,1}}=(\footnotesize{\underbrace{2,\ldots,2}_{k-1},x})\right)\nonumber\\
=&(1-p_{2L+1})\cdot\EE_{(\tiny{\underbrace{2,\ldots,2}_L,\underbrace{1,\ldots,1}_{k-L}})}(T_k-\tau_{k,1}~|~T_k>\tau_{k,1}).
\end{align}
Similarly, we can prove
\begin{equation}\label{e:pre2}
\EE_{(\tiny{\underbrace{2,\ldots,2}_L,\underbrace{1,\ldots,1}_{k-L}})}(\tau_{k,2}-\tau_{k,1}~|~T_k>\tau_{k,1})=\EE_{(\tiny{\underbrace{2,\ldots,2}_L,\underbrace{1,\ldots,1}_{k-L}})}(\tau_{k,2}-\tau_{k,1})=n^{k-1}
\end{equation}
and
\begin{equation}\label{e:pre3}
\EE_{(\tiny{\underbrace{2,\ldots,2}_L,\underbrace{1,\ldots,1}_{k-L}})}(\tau_{k,3}-\tau_{k,2}~|~T_k>\tau_{k,2})=\EE_{(\tiny{\underbrace{2,\ldots,2}_L,\underbrace{1,\ldots,1}_{k-L}})}(\tau_{k,3}-\tau_{k,2})=n^{k-1}
\end{equation}
together with the result of Lemma \ref{l:Etau_1}. Also, by Lemma \ref{l:same_dist}, we can show that
\begin{align}\label{e:pre4}
&\EE_{(\tiny{\underbrace{2,\ldots,2}_L,\underbrace{1,\ldots,1}_{k-L}})}(T_k-\tau_{k,2}~|~T_k>\tau_{k,1})\nonumber\\
=&\EE_{(\tiny{\underbrace{2,\ldots,2}_L,\underbrace{1,\ldots,1}_{k-L}})}(T_k-\tau_{k,2}~|~T_k>\tau_{k,2})\cdot\PP_{(\tiny{\underbrace{2,\ldots,2}_L,\underbrace{1,\ldots,1}_{k-L}})}(T_k>\tau_{k,2}~|~T_k>\tau_{k,1})\nonumber\\
=&\EE_{(\tiny{\underbrace{2,\ldots,2}_L,\underbrace{1,\ldots,1}_{k-L}})}(T_k-\tau_{k,2}~|~T_k>\tau_{k,2})\cdot\PP_{(\tiny{\underbrace{2,\ldots,2}_{k-1},1})}(T_k>\tau_{k,1})\nonumber\\
=&(1-p_{2k+1})\cdot\EE_{(\tiny{\underbrace{2,\ldots,2}_L,\underbrace{1,\ldots,1}_{k-L}})}(T_k-\tau_{k,2}~|~T_k>\tau_{k,2}).
\end{align}
By (\ref{e:pre1}) and (\ref{e:pre2}),
\begin{align}\label{e:fun1}
u&=\EE_{(\tiny{\underbrace{2,\ldots,2}_L,\underbrace{1,\ldots,1}_{k-L}})}(T_k-\tau_{k,1})=(1-p_{2L+1})\EE_{(\tiny{\underbrace{2,\ldots,2}_L,\underbrace{1,\ldots,1}_{k-L}})}(T_k-\tau_{k,1}~|~T_k>\tau_{k,1})\nonumber\\
&=(1-p_{2L+1})\left[\EE_{(\tiny{\underbrace{2,\ldots,2}_L,\underbrace{1,\ldots,1}_{k-L}})}(\tau_{k,2}-\tau_{k,1}~|~T_k>\tau_{k,1})+\EE_{(\tiny{\underbrace{2,\ldots,2}_L,\underbrace{1,\ldots,1}_{k-L}})}(T_k-\tau_{k,2}~|~T_k>\tau_{k,1})\right]\nonumber\\
&=(1-p_{2L+1})(n^{k-1}+v).
\end{align}
Similarly, by (\ref{e:pre3}) and (\ref{e:pre4}),
\begin{align}\label{e:fun2}
v&=\EE_{(\tiny{\underbrace{2,\ldots,2}_L,\underbrace{1,\ldots,1}_{k-L}})}(T_k-\tau_{k,2}~|~T_k>\tau_{k,1})=(1-p_{2k-1})\EE_{(\tiny{\underbrace{2,\ldots,2}_L,\underbrace{1,\ldots,1}_{k-L}})}(T_k-\tau_{k,2}~|~T_k>\tau_{k,2})\nonumber\\
&=(1-p_{2k-1})\left[\EE_{(\tiny{\underbrace{2,\ldots,2}_L,\underbrace{1,\ldots,1}_{k-L}})}(\tau_{k,3}-\tau_{k,2}~|~T_k>\tau_{k,2})+\EE_{(\tiny{\underbrace{2,\ldots,2}_L,\underbrace{1,\ldots,1}_{k-L}})}(T_k-\tau_{k,3}~|~T_k>\tau_{k,2})\right]\nonumber\\
&=(1-p_{2k-1})(n^{k-1}+v).
\end{align}
By (\ref{e:fun1}) and (\ref{e:fun2}), together with $(n-1)p_{2k-1}+p_{2k}=1$ which comes from the transitivity of $G_k$, we get
$$u=\frac{1-p_{2L+1}}{1-p_{2k}}\cdot(n-1)n^{k-1}.$$ Then together with Lemma \ref{l:Etau_1}, we get
\begin{align*}
s_{L,k}&=\EE_{(\tiny{\underbrace{2,\ldots,2}_L,\underbrace{1,\ldots,1}_{k-L}})}(T_k)\\
&=\EE_{(\tiny{\underbrace{2,\ldots,2}_L,\underbrace{1,\ldots,1}_{k-L}})}(\tau_{k,1})+\EE_{(\tiny{\underbrace{2,\ldots,2}_L,\underbrace{1,\ldots,1}_{k-L}})}(T_k-\tau_{k,1})\\
&=\frac{k}{k-1}s_{L,k-1}+\frac{1-p_{2L+1}}{1-p_{2k}}\cdot(n-1)n^{k-1}
\end{align*}
for $k\geq L+1$, as desired.\qed\\

Comparing (\ref{e:sk_induction2}) with (\ref{e:sk_induction}), we can see that we only need to calculate $p_{2k}$ and $p_{2L+1}$. The following lemma is the crucial step.
\begin{lemma}\label{l:difficult}
$p_1=p_{2k}$.
\end{lemma}

 \pf For any $m_1,m_2\in\{1,\ldots,2k\}$, let $$q_{m_1m_2}=\PP(X_1\in B_{m_2}~|~X_0\in B_{m_1}),$$ which equals to $\PP_x(X_1\in B_{m_2})$ for any $x\in B_{m_1}$ by the transitivity of $G_k$. Then for any $1\leq i\leq k$,
\begin{equation*}
\begin{cases}
q_{2i,2i-1}=\frac{1}{k};\\q_{2i,2i-2}=\frac{i-1}{k};\\q_{2i,2i+2}=\frac{k-i}{k}\cdot\frac{1}{n-1};\\q_{2i,2i}=\frac{k-i}{k}\cdot\frac{n-2}{n-1};\\
q_{2i-1,2i}=\frac{1}{k}\cdot\frac{1}{n-1};\\q_{2i-1,2i+1}=\frac{k-i}{k}\cdot\frac{1}{n-1};\\q_{2i-1,2i-3}=\frac{i-1}{k};\\q_{2i-1,2i-1}=\frac{k-i+1}{k}\cdot\frac{n-2}{n-1}.
\end{cases}
\end{equation*}
And $q_{m_1m_2}=0$ otherwise. From this, we first get
\begin{equation}\label{e:p2k}
p_{2k}=q_{2k,2k-2}p_{2k-2}=\frac{k-1}{k}p_{2k-2}.
\end{equation}
Next, $\{p_i:~i=1,\ldots,2k-2\}$ follows
\begin{equation}\label{e:equations}
\begin{cases}
p_1=\frac{n-2}{n-1}p_1+\frac{1}{k}\cdot\frac{1}{n-1}p_2+\frac{k-1}{k}\cdot\frac{1}{n-1}p_3;\\
p_2=\frac{1}{k}p_1+\frac{k-1}{k}\cdot\frac{n-2}{n-1}p_2+\frac{k-1}{k}\cdot\frac{1}{n-1}p_4;\\
p_3=\frac{1}{k}p_1+\frac{k-1}{k}\cdot\frac{n-2}{n-1}p_3+\frac{1}{k}\cdot\frac{1}{n-1}p_4+\frac{k-2}{k}\cdot\frac{1}{n-1}p_5;\\
~~~~~~~~~~~~~~~~~~~~~~~~~~~~~~~~~\vdots~~~~~~~~~~~~~~~~~~~~~~~~~~~~~~~~;\\
p_{2k-4}=\frac{k-3}{k}p_{2k-6}+\frac{1}{k}p_{2k-5}+\frac{2}{k}\cdot\frac{n-2}{n-1}p_{2k-4}+\frac{2}{k}\cdot\frac{1}{n-1}p_{2k-2};\\
p_{2k-3}=\frac{k-2}{k}p_{2k-5}+\frac{2}{k}\cdot\frac{n-2}{n-1}p_{2k-3}+\frac{2}{k}\cdot\frac{1}{n-1}p_{2k-2}.
\end{cases}
\end{equation}
Note that we do not need $p_{2k-1}$ in the last equation since if $\{X_t\}$ touches $B_{2k-1}$ before $B_{2k}$, then $T>\tau_1$. The first equation in (\ref{e:equations}) implies
\begin{equation}\label{e:eq1}
p_1=\frac{1}{k}p_2+\frac{k-1}{k}p_3.
\end{equation}
The second and third equations in (\ref{e:equations}) imply
\begin{equation}\label{e:eq2}
\left(1-\frac{k-1}{k}\cdot\frac{n-2}{n-1}\right)p_2=\frac{1}{k}p_1+\frac{k-1}{k}\cdot\frac{1}{n-1}p_4
\end{equation}
and
\begin{equation}\label{e:eq3}
\left(1-\frac{k-1}{k}\cdot\frac{n-2}{n-1}\right)p_3=\frac{1}{k}p_1+\frac{1}{k}\cdot\frac{1}{n-1}p_4+\frac{k-2}{k}\cdot\frac{1}{n-1}p_5.
\end{equation}
Putting (\ref{e:eq2}) and (\ref{e:eq3}) into (\ref{e:eq1}), we get
\begin{align*}
&k\left(1-\frac{k-1}{k}\cdot\frac{n-2}{n-1}\right)p_1\\
=&\left(\frac{1}{k}p_1+\frac{k-1}{k}\cdot\frac{1}{n-1}p_4\right)+(k-1)\left(\frac{1}{k}p_1+\frac{1}{k}\cdot\frac{1}{n-1}p_4+\frac{k-2}{k}\cdot\frac{1}{n-1}p_5\right).
\end{align*}
That is, $$p_1=\frac{2}{k}p_4+\frac{k-2}{k}p_5.$$ Similarly, we can inductively get
\begin{equation}\label{e:p1}
p_1=\frac{1}{k}p_2+\frac{k-1}{k}p_3=\frac{2}{k}p_4+\frac{k-2}{k}p_5=\ldots\ldots=\frac{k-1}{k}p_{2k-2}.
\end{equation}
The last equality is due to the absence of $p_{2k-1}$ as explained above. From (\ref{e:p2k}) and (\ref{e:p1}), we get $p_1=p_{2k}$, as desired.\qed\\

\noindent\textbf{Proof of Proposition \ref{p:sk_induction}.}~By the transitivity of $G_k$, we can get
\begin{equation}\label{e:extra}
(n-1)p_1+p_2=(n-1)p_3+p_4=(n-1)p_5+p_6=\ldots=(n-1)p_{2k-3}+p_{2k-2}=1.
\end{equation}
Together with (\ref{e:p1}), we obtain
$$p_{2j}=1-(n-1)p_{2j-1}=1-\frac{n-1}{k-j+1}[kp_1-(j-1)p_{2j-2}]$$ for $j=2,3,\ldots,k-1$. That is,
\begin{equation}\label{e:extraequation}
\begin{cases}
p_2=1-(n-1)p_1;\\p_4-\frac{n-1}{k-1}p_2=1-\frac{k(n-1)}{k-1}p_1;\\p_6-\frac{2(n-1)}{k-2}p_4=1-\frac{k(n-1)}{k-2}p_1;\\ ~~~~~~~~~~~~~~~~~~~\vdots~~~~~~~~~~~~~~~~~~~~~;\\p_{2k-2}-\frac{(k-2)(n-1)}{2}p_{2k-4}=1-\frac{k(n-1)}{2}p_1.
\end{cases}
\end{equation}
From (\ref{e:extraequation}),
\begin{align*}
p_{2k-2}&=\left[1-\frac{k(n-1)}{2}p_1\right]+\sum\limits_{i=1}^{k-2}\left(\prod\limits_{j=i}^{k-2}\frac{j(n-1)}{k-j}\right)\cdot\left[1-\frac{k(n-1)}{k-i+1}p_1\right]\\
&=\left[1-\frac{k(n-1)}{2}p_1\right]+\sum\limits_{i=1}^{k-2}\left[\frac{(k-2)!\cdot(n-1)^{k-i-1}}{(k-i)!(i-1)!}-\frac{k\cdot(k-2)!\cdot(n-1)^{k-i}}{(k-i+1)!(i-1)!}p_1\right]\\
&=\frac{1}{(k-1)(n-1)}\sum\limits_{i=1}^{k-1}\left[{k-1 \choose i-1}(n-1)^{k-i}-{k \choose i-1}(n-1)^{k-i+1}p_1\right]\\
&=\frac{1}{(k-1)(n-1)}[(n^{k-1}-1)-(n^k-1-k(n-1))p_1].
\end{align*}
Together with (\ref{e:p2k}),
\begin{align}\label{e:extra_last}
p_{2k}=\frac{k-1}{k}p_{2k-2}&=\frac{1}{k(n-1)}[(n^{k-1}-1)-(n^k-1-k(n-1))p_1]\nonumber\\
&=\frac{n^{k-1}-1}{k(n-1)}-\left[\frac{n^k-1}{k(n-1)}-1\right]p_1.
\end{align}
Lemma \ref{l:difficult} tells us that $p_1=p_{2k}$. Therefore, by (\ref{e:extra_last}),
\begin{equation}\label{e:p1p2k}
p_1=p_{2k}=\frac{n^{k-1}-1}{n^k-1}.
\end{equation}
Next, by (\ref{e:extraequation}),
\begin{align*}
p_{2L+2}&=\frac{1}{{k-1 \choose L}}\sum\limits_{i=1}^{L+1}\left[{k-1 \choose i-1}(n-1)^{L-i+1}-{k \choose i-1}(n-1)^{L-i+2}\cdot p_1\right]\\
&=1+\frac{1}{{k-1 \choose L}}\left[\sum\limits_{i=1}^L{k-1 \choose i-1}(n-1)^{L-i+1}-\sum\limits_{i=1}^{L+1}{k \choose i-1}(n-1)^{L-i+2}\cdot p_1\right].
\end{align*}
Together with (\ref{e:extra}),
\begin{align*}
p_{2L+1}&=\frac{1}{n-1}(1-p_{2L+2})\\
&=\frac{1}{{k-1 \choose L}}\left[\sum\limits_{i=1}^{L+1}{k \choose i-1}(n-1)^{L-i+1}\cdot p_1-\sum\limits_{i=1}^L{k-1 \choose i-1}(n-1)^{L-i}\right].
\end{align*}
Then together with (\ref{e:p1p2k}),
\begin{align}\label{e:term2}
&\frac{1-p_{2L+1}}{1-p_{2k}}\cdot(n-1)n^{k-1}=(n^k-1)(1-p_{2L+1})\nonumber\\
=&(n^k-1)+\frac{n^k-1}{{k-1 \choose L}}\sum\limits_{i=1}^L{k-1 \choose i-1}(n-1)^{L-i}-\frac{n^{k-1}-1}{{k-1 \choose L}}\sum\limits_{i=1}^{L+1}{k \choose i-1}(n-1)^{L-i+1}\nonumber\\
=&(n^k-1)+\frac{n^k-1}{{k-1 \choose L}}\sum\limits_{i=1}^L{k-1 \choose i-1}(n-1)^{L-i}-\frac{n^{k-1}-1}{{k-1 \choose L}}\sum\limits_{i=1}^{L+1}{k-1 \choose i-1}(n-1)^{L-i+1}\nonumber\\
&~~~~~~~~~~~~~~~~~~~~~~~~~~~~~~~~~~~~~~~~~~~~~-\frac{n^{k-1}-1}{{k-1 \choose L}}\sum\limits_{i=2}^{L+1}{k-1 \choose i-2}(n-1)^{L-i+1}\nonumber\\
=&(n^k-1)+\frac{n^k-1}{{k-1 \choose L}}\sum\limits_{i=1}^L{k-1 \choose i-1}(n-1)^{L-i}-\frac{n^{k-1}-1}{{k-1 \choose L}}\sum\limits_{i=1}^{L+1}{k-1 \choose i-1}(n-1)^{L-i+1}\nonumber\\
&~~~~~~~~~~~~~~~~~~~~~~~~~~~~~~~~~~~~~~~~~~~~~-\frac{n^{k-1}-1}{{k-1 \choose L}}\sum\limits_{i=1}^L{k-1 \choose i-1}(n-1)^{L-i}\nonumber\\
=&(n^k-1)-(n^{k-1}-1)\nonumber\\
&~+\left[\frac{1}{{k-1 \choose L}}\sum\limits_{i=1}^L{k-1 \choose i-1}(n-1)^{L-i}\right]\times[(n^k-1)-(n-1)(n^{k-1}-1)-(n^{k-1}-1)]\nonumber\\
=&n^{k-1}(n-1)+\frac{n-1}{{k-1 \choose L}}\sum\limits_{i=1}^L{k-1 \choose i-1}(n-1)^{L-i}.
\end{align}
The third equality is due to the combinatorial identity $\displaystyle{k \choose i-1}={k-1 \choose i-1}+{k-1 \choose i-2}$ for $i\geq2$. The fourth equality is due to the change of variable~(from $i$ to $i-1$) in the last term. Putting (\ref{e:term2}) into (\ref{e:sk_induction2}), we can get (\ref{e:sk_induction}), as desired.\qed

\bigskip

\noindent\textbf{Remark.}~When $L=0$, the recursive formula (\ref{e:sk_induction}) becomes
\begin{equation}\label{e:sk0_induction}
\begin{cases}
\displaystyle s_{0,k}=\frac{k}{k-1}s_{0,k-1}+(n-1)n^{k-1}, & k\geq 1;\\
s_{0,0}=0. &
\end{cases}
\end{equation}
It can be easily obtained from (\ref{e:sk0_induction}) that $$\EE_{(\tiny{\underbrace{1,\ldots,1}_M})}\left(T_{(\tiny{\underbrace{2,\ldots,2}_M})}\right)=s_{0,M}=M(n-1)\sum\limits_{k=1}^M\frac{n^{k-1}}{k}.$$ This is exactly the result in Corollary \ref{c:main}. Therefore, the recursive formula (\ref{e:sk0_induction}) provides a direct proof of Corollary \ref{c:main}. We can see that our present proof makes a good illustration for the recursive formula (\ref{e:sk0_induction}).

\section{Proof of Theorem \ref{t:main2}}

We can obtain from (\ref{e:sk_induction}) and (\ref{e:main}) that
\begin{align*}
&s_{L,M}=M(n-1)\sum\limits_{k=L+1}^M\frac{n^{k-1}}{k}+M(n-1)\sum\limits_{k=L+1}^M\sum\limits_{i=1}^L\frac{{k-1 \choose i-1}}{k{k-1 \choose L}}(n-1)^{L-i}\\
=&\sum\limits_{k=0}^{M-1}\frac{(n-1)^{k+1}}{{M-1 \choose k}}\sum\limits_{i=0}^k\frac{{M \choose i}}{(n-1)^i}-M(n-1)\sum\limits_{k=1}^L\frac{n^{k-1}}{k}+M(n-1)\sum\limits_{k=L+1}^M\sum\limits_{i=1}^L\frac{{k-1 \choose i-1}}{k{k-1 \choose L}}(n-1)^{L-i}.
\end{align*}
In order to prove $s_{L,M}=\displaystyle\sum\limits_{k=L}^{M-1}\frac{(n-1)^{k+1}}{{M-1 \choose k}}\sum\limits_{i=0}^k\frac{{M \choose i}}{(n-1)^i}$, we need to prove $$\sum\limits_{k=0}^{L-1}\frac{(n-1)^{k+1}}{{M-1 \choose k}}\sum\limits_{i=0}^k\frac{{M \choose i}}{(n-1)^i}=M(n-1)\sum\limits_{k=1}^L\frac{n^{k-1}}{k}-M(n-1)\sum\limits_{k=L+1}^M\sum\limits_{i=1}^L\frac{{k-1 \choose i-1}}{k{k-1 \choose L}}(n-1)^{L-i},$$ that is,
\begin{equation}\label{e:goal1}
\sum\limits_{k=0}^{L-1}\sum\limits_{i=0}^k(n-1)^{k-i}\cdot\frac{{M \choose i}}{M{M-1 \choose k}}+\sum\limits_{k=L+1}^M\sum\limits_{i=1}^L\frac{{k-1 \choose i-1}}{k{k-1 \choose L}}(n-1)^{L-i}=\sum\limits_{k=1}^L\frac{n^{k-1}}{k}.
\end{equation}
Note that $$\text{(LHS) of}~(\ref{e:goal1})=\sum\limits_{i=0}^{L-1}(n-1)^i\cdot\sum\limits_{k=i}^{L-1}\frac{{M \choose k-i}}{M{M-1 \choose k}}+\sum\limits_{i=0}^{L-1}(n-1)^i\cdot\sum\limits_{k=L+1}^M\frac{{k-1 \choose L-i-1}}{k{k-1 \choose L}},$$ and by (\ref{e:main})~(change $M$ by $L$), \begin{align*}
\text{(RHS) of}~(\ref{e:goal1})&=\sum\limits_{k=0}^{L-1}\frac{1}{L{L-1 \choose k}}\cdot\sum\limits_{i=0}^k{L \choose i}(n-1)^{k-i}\\
&=\sum\limits_{k=0}^{L-1}\frac{1}{L{L-1 \choose k}}\cdot\sum\limits_{i=0}^k{L \choose k-i}(n-1)^i=\sum\limits_{i=0}^{L-1}(n-1)^i\cdot\sum\limits_{k=i}^{L-1}\frac{{L \choose k-i}}{L{L-1 \choose k}}.
\end{align*}
So to prove (\ref{e:goal1}), it suffices to prove that for any $i\in\{0,\ldots,L-1\}$,
\begin{equation}\label{e:goal2}
\sum\limits_{k=i}^{L-1}\frac{{M \choose k-i}}{M{M-1 \choose k}}+\sum\limits_{k=L+1}^M\frac{{k-1 \choose L-i-1}}{k{k-1 \choose L}}=\sum\limits_{k=i}^{L-1}\frac{{L \choose k-i}}{L{L-1 \choose k}}.
\end{equation}
Note that
\begin{align*}
\sum\limits_{k=i}^{L-1}\frac{{L \choose k-i}}{L{L-1 \choose k}}-\sum\limits_{k=i}^{L-1}\frac{{M \choose k-i}}{M{M-1 \choose k}}&=\sum\limits_{k=i}^{L-1}\frac{{k \choose i}}{(i+1){L-k+i \choose i+1}}-\sum\limits_{k=i}^{L-1}\frac{{k \choose i}}{(i+1){M-k+i \choose i+1}}\\
&=\sum\limits_{k=i}^{L-1}\frac{{k\choose i}}{i+1}\sum\limits_{j=L+1}^M\left[\frac{1}{{j-k+i-1 \choose i+1}}-\frac{1}{{j-k+i \choose i+1}}\right]\\
&=\sum\limits_{j=L+1}^M\sum\limits_{k=i}^{L-1}\frac{(i+1)k!(j-k-2)!}{(k-i)!(j-k+i)!},
\end{align*}
and $$\sum\limits_{k=L+1}^M\frac{{k-1 \choose L-i-1}}{k{k-1 \choose L}}=\sum\limits_{j=L+1}^M\frac{L!(j-L-1)!}{j(L-i-1)!(j-L+i)!}.$$
Therefore, to prove (\ref{e:goal2}), it suffices to prove that for any $j\in\{L-1,\ldots,M\}$,
$$\sum\limits_{k=i}^{L-1}\frac{(i+1)k!(j-k-2)!}{(k-i)!(j-k+i)!}=\frac{L!(j-L-1)!}{j(L-i-1)!(j-L+i)!},$$
which suffices to prove that for any $L\geq i+1$, $$\frac{(i+1)L!(j-L-2)!}{(L-i)!(j-L+i)!}=\frac{(L+1)!(j-L-2)!}{j(L-i)!(j-L+i-1)!}-\frac{L!(j-L-1)!}{j(L-i-1)!(j-L+i)!},$$
which is easy to check. Therefore, (\ref{e:goal2}) holds, and furthermore, (\ref{e:goal1}) holds. The proof of Theorem \ref{t:main2} is now complete.\qed

\section{Concluding remarks}

In this paper, we use a new, purely probabilistic method--the method of stopping times, to compute the expected hitting time of moving from a given configuration to a different one for the multiple--urn Ehrenfest model. As a special case, we compute the expected hitting time when all $M$ balls are placed in a specific urn given that initially all $M$ balls are in another urn. This extends the result in Chen et al.~\cite{Chen-Goldstein-Lathrop-Nelsen2017}.

For some technical reasons, we only consider the ``unbiased'' case, that is, the ball is chosen randomly and put in other urns randomly with equal probabilities. It will be interesting to consider the case of biased or preferential probabilities in the future.

Furthermore, we know that Lemma \ref{l:difficult} is a crucial step in the proof. An interesting problem is: can we give a bijective proof to illustrate this ``simple'' equality immediately?

\bigskip

\noindent\textbf{Acknowledgment.}\quad We would like to thank two anonymous referees for their useful suggestions,
which are of great help for improving the manuscript. The research of Qiang Yao was partially supported by the Natural Science Foundation of China~(No.11671145), the program of China Scholarships Council~(No.201806145024), and the 111 Project~(B14019).


\bigskip


\begin{thebibliography}{99}

\bibitem{Aloisi-Nali2018}Aloisi, A.M. and Nali, P.F. (2018). Marbles and bottles--or boxes illustrate irreversibility and recurrence. \emph{Physics Education(India)}, \textbf{34}, 1--18.

\bibitem{Blom1989}Blom, G. (1989). Mean transition times for the Ehrenfest urn model. \emph{Adv. Appl. Prob.} \textbf{21}, 479--480.

\bibitem{Chen-Goldstein-Lathrop-Nelsen2017} Chen, Y-P., Goldstein, I.H., Lathrop, E.D. \& Nelsen, R.B. (2017). Computing an expected hitting time for the 3-urn Ehrenfest model via electric networks. {\em Stat. Prob. Letters} \textbf{127}, 42--48.

\bibitem{Doyle-Snell1984}Doyle, P.G. and Snell, E.J. (1984). Random Walks and Electric Networks. \emph{Carus Math.
Monographs}, \textbf{22}, Math. Assoc. Amer., Washington, D.C.


\bibitem{Ehrenfest-Ehrenfest1907}Ehrenfest, P. and Ehrenfest, T. (1907). \"Uber zwei bekannte Einw\"ande gegen das Boltzmannsche $H$--Theorem.
\emph{Physikalische Zeitschrift} \textbf{8}, 311--314.

\bibitem{Lathrop-Goldstein-Chen2016}Lathrop, E.D., Goldstein, I.H. and Chen, Y-P. (2016). A note on a generalized Ehrenfest urn model: another look at the mean transition times. \emph{J. Appl. Prob.} \textbf{53}, 630--632.

\bibitem{Lyons-Peres2017}Lyons, R. and Peres, Y. (2017). Probability on Trees and Networks. \emph{Cambridge Series in Statistical and Probabilistic Mathematics}, \textbf{42}, Cambridge University Press, New York.

\bibitem{Meerson-Zilber2018}Meerson, B. and Zilber, B. (2018). Large deviations of a long-time average in the Ehrenfest urn model. \emph{Journal of Statistical Mechanics: Theory and Experiment}, \textbf{2018(5)}, 053202.

\bibitem{Palacios1994}Palacios, J.L. (1994). Another look at the Ehrenfest urn via electric networks. \emph{Adv. Appl. Prob.} \textbf{26}, 820--824.

\bibitem{Xin-Zhao-Yao-Cui2020}Xin, C., Zhao, M., Yao, Q. and Cui, E. (2020). On the distribution of the hitting time for the $N$--urn Ehrenfest model. {\em Stat. Prob. Letters} \textbf{157}, 108625.

\bibitem{Xue2021}Xue, X. (2021). Hydrodynamics of the generalized $N$--urn Ehrenfest model. \emph{Preprint}, current version available at \emph{https://arxiv.org/pdf/2010.08726.pdf}.

\end{thebibliography}
\end{document}